\documentclass[nopreprintline,preprint]{elsarticle}

\journal{Journal of Computational Physics}

\usepackage{lineno,hyperref}
\modulolinenumbers[5]

\usepackage[utf8]{inputenc}

\usepackage{amsmath,amsfonts,amsthm,amssymb,bm} 

\usepackage{xcolor}

\usepackage{siunitx} 
\usepackage{booktabs}

\def\appendixname{}
\usepackage{cleveref}
\crefname{lstlisting}{listing}{listings}
\Crefname{lstlisting}{Listing}{Listings}
\crefname{figure}{figure}{figures}
\crefformat{equation}{(#2#1#3)}
\crefrangeformat{equation}{(#3#1#4) to~(#5#2#6)}
\crefmultiformat{equation}{(#2#1#3)}%
{ and~(#2#1#3)}{, (#2#1#3)}{ and~(#2#1#3)}
\usepackage{listings}

\usepackage[makeroom]{cancel}

\usepackage{subcaption}

\usepackage{pgfplots}
\usetikzlibrary{calc}
\usepackage{tikz-3dplot}
\pgfplotsset{compat=1.17}

\renewcommand\vec{\mathbf}

\newcommand{\avg}[1]{\{#1\}}
\newcommand{\jump}[1]{[\![#1]\!]}

\newcommand*\chem[1]{\ensuremath{\mathrm{#1}}}

\newcommand\numberthis{\addtocounter{equation}{1}\tag{\theequation}}

\makeatletter
\DeclareRobustCommand{\pder}[1]{%
  \@ifnextchar\bgroup{\@pder{#1}}{\@pder{}{#1}}}
\newcommand{\@pder}[2]{\frac{\partial#1}{\partial#2}}
\makeatother

\makeatletter
\DeclareRobustCommand{\der}[1]{%
  \@ifnextchar\bgroup{\@der{#1}}{\@der{}{#1}}}
\newcommand{\@der}[2]{\frac{d#1}{d#2}}
\makeatother

\newcommand{\RR}{\mathbb{R}}

\newcommand{\mb}[1]{\mathbf{#1}}

\newcommand*\diff{\mathop{}\!\mathrm{d}}

\definecolor{ao(english)}{rgb}{0.0, 0.5, 0.0}
\definecolor{electricpurple}{rgb}{0.75, 0.0, 1.0}
\newcommand{\rone}[1]{#1}
\newcommand{\rtwo}[1]{#1}
\newcommand{\rme}[1]{#1}

\begin{document}

\begin{frontmatter}
\title{A scalable DG solver for the electroneutral Nernst-Planck equations}

\author{Thomas Roy\corref{cor1}}
\ead{roy27@llnl.gov}
\author{Julian Andrej}
\ead{andrej1@llnl.gov}
\author{Victor A. Beck}
\ead{beck33@llnl.gov}
\address{Lawrence Livermore National Laboratory, Livermore, 94550, CA, USA}
\cortext[cor1]{Corresponding author}
\begin{abstract}
    The robust, scalable simulation of flowing electrochemical systems is increasingly important due to the synergy between intermittent renewable energy and electrochemical technologies such as energy storage and chemical manufacturing.
    The high P\'eclet regime of many such applications prevents the use of off-the-shelf discretization methods.
    In this work, we present a high-order Discontinuous Galerkin scheme for the electroneutral Nernst-Planck equations.
    The chosen charge conservation formulation allows for the specific treatment of the different physics: upwinding for advection and migration, and interior penalty for diffusion of ionic species as well the electric potential.
    Similarly, the formulation enables different treatments in the preconditioner: AMG for the potential blocks and ILU-based methods for the advection-dominated concentration blocks.
    We evaluate the convergence rate of the discretization scheme through numerical tests.
    Strong scaling results for two preconditioning approaches are shown for a large 3D flow-plate reactor example.
\end{abstract}

\end{frontmatter}



\section{Introduction}

With the steady transition to renewable energy comes a surge of intermittent electricity \cite{Chu_2012,Chu_2016}.
Fortunately, various electrochemical technologies such as energy storage and electrochemical manufacturing can capitalize on this otherwise wasted energy \cite{Gur_2018,Ager_2018}.
To handle the rise of renewable energy, the development and scaling-up of these technologies is crucial.
Naturally, this brings a need for simulation to aid design of electrochemical devices at the industrial scale.

Such simulations require the solution of highly coupled systems of partial differential equations (PDEs).
An important model in the presence of fluid flow is the Nernst-Planck (NP) equation, which describes the transport of an ionic species by advection, diffusion and electromigration \cite{newman2012electrochemical}.
For multi-ion transport, the system consisting of an NP equation for each species can be closed by a Poisson equation for the electric potential.
However, the Poisson equation operates on a much smaller lengthscale and is thus challenging to solve computationally \cite{dickinson2011electroneutrality}.
Due to the scale separation, it is commonly replaced by the electroneutrality condition: an algebraic equation for the ionic concentrations.

For the electroneutral system, arriving at a steady-state solution is a challenging numerical
task and therefore the discretization scheme has to be chosen with
care. In particular, specific treatment is required for the typically advection-dominated nature of the flow.
The considered physical equations have been successfully discretized with
several methods in the past. The authors in~\cite{bortels1996multi} use a
multi-dimensional upwinding method similar to a finite volume approach. A
continuous Galerkin finite element method has to employ a stabilization scheme
because of the advection-dominated system, which was shown
in~\cite{bauer2012coupled} using a Variational Multiscale Method. There have
been attempts at one-dimensional problems using the Discontinuous Galerkin (DG)
Method~\cite{liu2017free,sun2018discontinuous}, which showed promising results and motivated this research. To the
best of the authors' knowledge, none of the previously mentioned references
looked into the benefits of high-order DG spatial discretizations, which are
frequently used in advection-dominated transport and compressible flow
applications~\cite{hartmann2014higher}.

After discretization and linearization, large, sparse linearized systems must be solved.
Iterative methods using appropriate preconditioning techniques \cite{wathen2015preconditioning} are key to obtaining a scalable solver.
Designing preconditioners for problems with several, simultaneous, coupled physical phenomena is especially challenging since different parts of the problem may require different treatments.
In the case of electroneutral ion transport, the equations are elliptic with respect to the electrical potential but hyperbolic with respect to the concentration variables.
Scalable solvers for this problem are not available in commercial tools, where directs solver are used.
This is especially restrictive when aiming for industrial scale 3D simulations,
as well as the simulation of porous electrodes at the resolved pore-scale.

In this paper, we present a high-order DG spatial discretization and scalable preconditioning approaches for the electroneutral Nernst-Planck equations.
The electroneutral system is arranged in a charge-conservation formulation that is necessary for the proper penalization of the potential terms in the DG scheme.
Additionally, this formulation also allows for block preconditioning with specific treatment of the potential block and concentration blocks.

In \cref{sec:gov},
the equations for the relevant electrochemical physics in an electroneutral
setting are described. We explain how to reorder the equations and exploit the
electroneutrality condition to eliminate a species variable and construct the charge conservation
equation. Next, in \cref{sec:dg}, we detail the spatial discretization
process using the DG finite element method including
necessary treatment for the hyperbolic and elliptic parts of the resulting system.
Then, in \cref{sec:pc}, we outline our approach to numerically precondition and solve the
resulting discrete nonlinear problem. Finally, in \cref{sec:results}, we perform
a convergence test for our discrete scheme and show scaling results for a large
three-dimensional, industry-relevant problem.

\section{Governing Equations}
\label{sec:gov}

We consider the equations for electroneutral multi-ion transport. First, we
describe the standard Electroneutral-Nernst-Planck model, as well as
appropriate boundary conditions for a flow-plate reactor. Then, we discuss
alternative formulations including the Charge-Conservation-Nernst-Planck model,
which will be used for our Discontinuous Galerkin scheme in \cref{sec:dg}. A nondimensionalization
of that model is also given.

\subsection{Electroneutral Nernst-Planck equation}
\label{subsec:gov_npe}

We consider the flow of a dilute electrolyte solution containing $m\geq 2$
chemical species in a domain $\Omega \in \RR^d$, $d=1,2$, or 3. In the presence
of an electric field, the mass transport is described by the Nernst-Planck
equation \cite{newman2012electrochemical}.
For simplicity, we do not consider the transient case, although our approach would still be applicable.

At steady state, the conservation of mass of each species $k=1,\dots,m$ is given by
\begin{equation}\label{eq:NP}
    -\nabla \cdot (D_k \nabla c_k) + \nabla \cdot (c_k\mb{u} ) - \nabla \cdot (z_k
    \mu_k F c_k \nabla \Phi) - R_k = 0 \quad \text{in}\quad\Omega,
\end{equation}
where $c_k$ is the molar concentration of species $k$, $D_k$ is the molecular
diffusion coefficient, $\mb{u}$ is the velocity of the electrolyte, $z_k$ is the
valence (charge number), $\mu_k$ is the mobility constant, $F$ is
the Faraday constant, $\Phi$ is the electric potential in the electrolyte,
and $R_k$ is a reaction term. The first term in \cref{eq:NP} describes molecular diffusion, the second,
advection, and the third, electromigration.

The mobility constant is given by the Nernst-Einstein relation \cite{newman2012electrochemical}
\begin{equation}\label{eq:NE}
    \mu_k = \frac{D_k}{RT},
\end{equation}
where $R$ is the universal gas constant, and $T$ is the thermodynamic temperature.

The system currently consists of $m$ equations \cref{eq:NP}, and $m+1$
unknowns: the $m$ concentrations $c_k$ and the potential $\Phi$. An additional
equation is required to complete the system. A commonly-used relation is the
electroneutrality approximation
\begin{equation}\label{eq:electroneutrality}
    \sum_{k=1}^m z_k c_k = 0 \quad \text{in}\quad \Omega,
\end{equation}
which states that there is no charge separation and that neutrality is
maintained everywhere. This approximation is obtained from the Poisson equation
and is appropriate when considering length scales greater than nanometers
\cite{dickinson2011electroneutrality}.

In \cite{bauer2012coupled}, the system consisting of \cref{eq:NP,eq:electroneutrality} is called the Electroneutrality-Nernst-Planck (ENP) model.

For convenience, we denote the advection-migration transport as $\mb{q}_k = \mb{u} - z_k
\mu_k F\nabla \Phi$ so that \cref{eq:NP} can be written as
\begin{equation}\label{eq:NP2}
    -\nabla \cdot (D_k \nabla c_k) + \nabla \cdot (c_k\mb{q}_k ) - R_k = 0 \quad \text{in}\quad\Omega.
\end{equation}
We now provide boundary conditions that are appropriate for a flow plate reactor.
The boundary $\partial\Omega$ is partitioned with respect to the direction of
the velocity at the boundary into the inlet boundary $\Gamma_\mathrm{in}$, the
outlet boundary $\Gamma_\mathrm{out}$, the wall boundary $\Gamma_\mathrm{wall}$
and the electrode boundary $\Gamma_\mathrm{elec}$. We have
\begin{align*}
    \Gamma_\mathrm{in} &= \left\{ \mb x \in \partial \Omega \;\vert\; \mb u \cdot \mb
    n < 0 \right\}, \\
    \Gamma_\mathrm{out} &= \left\{ \mb x \in \partial \Omega \;\vert\; \mb u \cdot \mb
    n > 0 \right\}, \\
    \Gamma_\mathrm{wall} \cup \Gamma_\mathrm{elec} &= \left\{ \mb x \in \partial \Omega \;\vert\; \mb u \cdot \mb
    n = 0 \right\},
\end{align*}
where $\mb n$ is the outward normal,
so that $\partial \Omega = \Gamma_\mathrm{in}\cup \Gamma_\mathrm{out} \cup
\Gamma_\mathrm{wall} \cup \Gamma_\mathrm{elec}$. Each boundary is such that
$\Gamma_i \cap \Gamma_j = \varnothing$. The location of the wall and electrode
boundaries are determined beforehand.

The boundary conditions are given by
\begin{equation}\label{eq:inlet}
    (c_k \mb q_k- D_k \nabla c_k)\cdot \mb n = c_k^\mathrm{in} \mb u \cdot \mb
    n \quad \text{on}\quad \Gamma_\mathrm{in},
\end{equation}
\begin{equation}\label{eq:outlet}
    (D_k \nabla c_k + z_k \mu_k F c_k \nabla \Phi)\cdot \mb n = 0
     \quad \text{on} \quad \Gamma_\mathrm{out},
\end{equation}
\begin{equation}\label{eq:electrode}
    (D_k \nabla c_k + z_k \mu_k F c_k \nabla \Phi)\cdot \mb n =  g_{k}\quad
    \text{on} \quad \Gamma_\mathrm{elec},
\end{equation}
\begin{equation}\label{eq:wall}
    (D_k \nabla c_k + z_k \mu_k F c_k \nabla \Phi)\cdot \mb n = 0 \quad
    \text{on} \quad \Gamma_\mathrm{wall},
\end{equation}
where $g_k$ represents electrochemical reactions at the electrode surface.

In our numerical experiments, we will consider the specific case of a simple $n$-electron transfer redox reaction $Ox + n e^- \rightarrow Red$, where the current density $J_n$ at the electrodes is given by the
Butler-Volmer equation
\begin{equation}\label{eq:bv}
    J_n = J_0  \left[\left(\frac{c_r}{c_r^*}\right)^\gamma
        \exp \left( \alpha_1 \frac{nF}{RT} (\Phi_\mathrm{app} -\Phi)\right)
    -\left(\frac{c_o}{c_o^*}\right)^\gamma\exp \left( -\alpha_2 \frac{nF}{RT} (\Phi_\mathrm{app} -\Phi)\right)\right],
\end{equation}
where $J_0$ is the exchange current density, $c_o$ and $c_r$ are the concentrations of the oxidant $Ox$ and reductant $Red$, respectively, $c_o^*$ and $c_r^*$ are reference concentrations, $\alpha_1$ and $\alpha_2$ are the anodic and cathodic charge transfer coefficients, respectively, $n$ is the number of electrons in the reaction, and $\Phi_\mathrm{app}$ is the applied potential.
The reaction rate in \cref{eq:electrode} is given by $g_k = \frac{J_n}{nF}$ if $k$ is the oxidant, $g_k=-\frac{J_n}{nF}$ if it is the reductant, and $g_k=0$ otherwise.

\subsection{Alternative formulations}
\label{subsec:gov_alt}
Alternative formulations of the ENP system \cref{eq:NP,eq:electroneutrality}
have different advantages with respect to the numerical solution
\cite{bauer2012coupled}.

In our case, we are interested in the Discontinuous Galerkin (DG) finite element
method \cite{riviere2008discontinuous}. In particular, we will consider upwinding techniques for the
advection and migration terms. The current formulation, ENP, is problematic with
respect to the stabilization and penalization of the potential $\Phi$. Indeed,
the mass conservation equations \cref{eq:NP} will be treated as transport
problems for the concentrations, using an upwinding scheme, leaving $\Phi$
untreated. In this electroneutral Nernst-Planck system, the electroneutrality
approximation replaces the Poisson equation, where stabilization and
penalization appropriate for an elliptic problem could have been used. We need
an equation that plays this role for $\Phi$.

We formulate a charge conservation equation by summing over all conservation of
mass equations \cref{eq:NP} weighting each by $z_k F$. The advection terms
cancel due to electroneutrality \cref{eq:electroneutrality}, and we obtain
\begin{equation}\label{eq:charge}
    -F \nabla \cdot \left[ \sum_{k=1}^m z_k D_k \nabla c_k\right] - \nabla \cdot
    \left[\left(\sum_{k=1}^m z_k^2 \mu_k F^2 c_k\right) \nabla \Phi\right] -
    \sum_{k=1}^m z_k F R_k = 0 \quad \text{in}\quad \Omega .
\end{equation}
However, we note that the system consisting of \cref{eq:NP} and
\cref{eq:charge} is linearly dependent. To circumvent this, one of the
concentration unknowns can be eliminated using the electroneutrality condition.
Without loss of generality, we eliminate the concentration $c_m$ using the
substitution
\begin{equation}\label{eq:substitution}
    c_m = -\frac{1}{z_m} \sum_{k=1}^{m-1} z_k c_k,
\end{equation}
assuming that $z_m \neq 0$. Combining \cref{eq:charge,eq:substitution}, we obtain
\begin{multline*}
    -F \nabla \cdot \left[ \sum_{k=1}^{m-1} z_k (D_k-D_m) \nabla c_k\right] - \nabla \cdot
    \left[\left(\sum_{k=1}^{m-1} z_k (z_k\mu_k-z_m\mu_m) F^2 c_k\right) \nabla
    \Phi\right] \\- \sum_{k=1}^{m-1} z_k F R_k = 0 \quad \text{in}\quad \Omega .\numberthis\label{eq:charge2}
\end{multline*}
This conservation of charge equation replaces the conservation of mass of
species $m$ in the system. The system now consists of $m$ equations:
conservation of mass equations \cref{eq:NP} for $k=1,\dots,m-1$, and a
conservation of charge equation \cref{eq:charge2}, and $m$ unknowns:
concentration unknowns $c_k$ for $k=1,\dots,m-1$, and the electric potential
$\Phi$. As in \cite{bauer2012coupled}, we call this formulation the
Charge-conservation-Nernst-Planck (CNP) model.

In terms of the DG scheme, stabilization and penalization of $c_m$ is no longer
needed since it is no longer an unknown of the system. We can instead treat
\cref{eq:charge2} as an elliptic equation for $\Phi$ in terms of the DG
stabilization and penalization.

\subsection{Nondimensionalization}
\label{sec:nondim}
Here we perform a nondimensionalization of the CNP model consisting of
\cref{eq:NP,eq:charge} and boundary conditions
\cref{eq:inlet,eq:outlet,eq:electrode,eq:wall}. Note that this scaling process
generally helps numerical algorithms. Indeed, for the considered test cases,
this greatly accelerates the convergence of the iterative methods.

Let the following nondimensional variables:
\begin{align}\label{eq:nondims}
    \hat{\Phi} = \frac{F}{R T} \Phi ,\quad
    \hat{c}_k = \frac{c_k}{c_k^\mathrm{in}},\quad
    \hat{\mb x} = \frac{1}{L}\mb x, \quad
    \hat{\mb u} = \frac{1}{u_\mathrm{avg}}\mb u,
\end{align}
and the differential operator with respect $\hat{\mb x}$
\begin{equation}\label{eq:nablahat}
    \hat\nabla = L \nabla.
\end{equation}
Substituting \cref{eq:nondims,eq:nablahat} into \cref{eq:NP}, we obtain
\begin{equation}\label{eq:NPND}
    -\hat\nabla \cdot (\hat D_k \hat\nabla \hat c_k) + \hat\nabla \cdot (\hat c_k \hat{\mb u} ) - \hat\nabla \cdot (z_k
    \hat D_k \hat c_k \hat\nabla \hat\Phi) - \frac{L}{c_k^\mathrm{in}u_\mathrm{avg}}R_k = 0 \quad \text{in}\quad\Omega,
\end{equation}
where $\hat D_k = \frac{D_k}{L u_\mathrm{avg}}$ is the inverse P\'eclet number.
We note that the electroneutrality condition \cref{eq:electroneutrality} is now
\begin{equation}
    \sum_{k=1}^m z_k c_k^\mathrm{in} \hat{c}_k = 0,
\end{equation}
so that $\hat{c}_m$ is eliminated using
\begin{equation}
    \hat{c}_m = - \frac{1}{z_m c_m^\mathrm{in}} \sum_{k=1}^{m-1} z_k c_k^\mathrm{in} \hat{c}_k.
\end{equation}
In order to nondimensionalize the charge conservation equation \cref{eq:charge},
we introduce the reference concentration $c_\mathrm{ref}$ to obtain
\begin{multline}\label{eq:chargeND}
    -\hat\nabla \cdot \left[ \sum_{k=1}^m z_k \hat D_k \frac{c_k^\mathrm{in}}{c_\mathrm{ref}}\hat\nabla \hat c_k\right] - \hat\nabla \cdot
    \left[\left(\sum_{k=1}^m z_k^2 \hat D_k \frac{c_k^\mathrm{in}}{c_\mathrm{ref}}\hat c_k\right) \hat\nabla \hat\Phi\right] \\
    -\sum_{k=1}^m \frac{L}{c_\mathrm{ref} u_\mathrm{avg}} z_k R_k = 0 \quad \text{in}\quad \Omega .
\end{multline}
The boundary conditions \cref{eq:inlet,eq:outlet,eq:electrode,eq:wall} become
\begin{equation}
    (\hat c_k \hat{\mb q}_k- \hat D_k \hat \nabla c_k)\cdot \mb n = \hat{\mb u} \cdot \mb
    n \quad \text{on}\quad \Gamma_\mathrm{in},
\end{equation}
\begin{equation}
    (\hat D_k \hat\nabla \hat c_k + z_k \hat D_k \hat c_k \hat \nabla \hat \Phi)\cdot \mb n = 0
     \quad \text{on} \quad \Gamma_\mathrm{out},
\end{equation}
\begin{equation}\label{eq:electrodeND}
    (\hat D_k \hat \nabla \hat{c}_k + z_k \hat D_k \hat c_k \hat\nabla
    \hat\Phi)\cdot \mb n =  \frac{1}{u_\mathrm{avg} c_k^\mathrm{in}}g_{k}\quad
    \text{on} \quad \Gamma_\mathrm{elec},
\end{equation}
\begin{equation}
    (\hat D_k \hat \nabla \hat c_k + z_k \hat D_k \hat c_k \hat \nabla \hat \Phi)\cdot \mb n = 0 \quad
    \text{on} \quad \Gamma_\mathrm{wall}.
\end{equation}
In the specific case of the Butler-Volmer equation \cref{eq:bv}, if $k$ is an oxidant, \cref{eq:electrodeND} becomes
\begin{multline}
    (\hat D_k \hat \nabla \hat c_k + z_k \hat D_k \hat c_k \hat\nabla
    \hat\Phi)\cdot \mb n =   \hat J_0 \left[\left(\hat c_r\frac{c_r^\mathrm{in}}{c_r^*}\frac{c_k^*}{c_k^\mathrm{in}}\right)^\gamma\exp \left(\hat{\Phi}_\mathrm{app} - \hat\Phi\right)
    \right. \\ \left.
    - \hat c_k^\gamma\exp \left( -(\hat{\Phi}_\mathrm{app} - \hat\Phi)\right)\right],
\end{multline}
where $\hat J_0 = \dfrac{J_0 (c_k^\mathrm{in})^{\gamma-1}}{u_\mathrm{avg} F(c_k^*)^\gamma }$ and $\hat{\Phi}_\mathrm{app} = \dfrac{F}{RT} \Phi_\mathrm{app}$.


\section{Discontinuous Galerkin method}
\label{sec:dg}
In this section, we introduce a Discontinuous Galerkin (DG) scheme for the CNP
model \cref{eq:NP2,eq:charge2}. After reviewing preliminary DG concepts, we
provide DG schemes for a hyperbolic and an elliptic equation. The resulting
numerical fluxes are then used for the full system \cref{eq:NP2,eq:charge2}. For
generality, we provide the scheme for the dimensional equations instead of the
nondimensionalized version from \cref{sec:nondim}.

\subsection{Preliminaries}
We begin by introducing the appropriate functional notation. Let $\mathcal{T}_h$
be a subdivision of $\Omega$ with element $K \in
\mathcal{T}_h$. The broken Sobolev spaces are defined as
\begin{align}
    H^l (\mathcal{T}_h) = \{v \in L^2(\Omega) : v|_{K} \in H^l (K) \}.
\end{align}
Let $\Gamma= \cup_{K\in\mathcal{T}_h} \partial K$ denote the set of element
boundaries. The restrictions of functions in $H^1(\mathcal{T}_h)$ to $\Gamma$
are called \emph{traces} and belong to the function space $T(\Gamma)=
\Pi_{K\in\mathcal{T}_h} L^2(\partial K)$.

Denote the space of polynomials of degree $p$ on each element
$\kappa$ with $Q^p$. A proper finite element space is defined as
\begin{align}
    V^p_h = \{v_h \in L^2(\Omega) : v_h|_{K} \in Q^p \}.
\end{align}

Next we introduce some trace operators needed for the numerical fluxes. For $u \in T(\Gamma)$, we define the average $\avg{u}$ and jump $\jump{u}$ on interior faces as
\begin{align}
    \avg{u} = \frac{1}{2} (u^- + u^+), \quad
    \jump{u} = u^- \vec{n}^- + u^+ \vec{n}^+ \quad \text{ on } e \in \mathcal{E}^\circ,
\end{align}
and for $\vec{q} \in [T(\Gamma)]^d$,
\begin{align}
    \avg{\vec{q}} = \frac{1}{2} (\vec{q}^- + \vec{q}^+), \quad
    \jump{u} = \vec{q}^- \cdot \vec{n}^- + \vec{q}^+ \cdot \vec{n}^+ \quad \text{ on } e \in \mathcal{E}^\circ.
\end{align}
Jumps on exterior faces are defined as
\begin{align}
    \jump{u} = u \vec{n}, \quad \avg{\vec{u}} = \vec{u} \quad \text{ on } e \in \mathcal{E}^\Gamma.
\end{align}

\subsection{Applied to the CNP model}
Following a combination of~\cite{brezzi2004discontinuous,cockburn2001runge,arnold2002unified}, we derive a
Discontinuous Galerkin formulation for the CNP model. In this section, we only
consider zero-Dirichlet boundary conditions for simplicity.  The full scheme
with boundary conditions \cref{eq:inlet,eq:outlet,eq:electrode,eq:wall} is given
in \cref{subsec:dg_scheme}.

Analyzing \cref{eq:NP2} we can see that it consists of a diffusion,
advection-migration, and reaction term. Each term has to be treated carefully in
order to achieve a convergent DG discretization.

We multiply~\cref{eq:NP2} by a test function $v$ and integrate over the domain
$\Omega$.
\begin{align}
    -\int_\Omega \nabla \cdot (D_k \nabla c_k) v \diff x
    + \int_\Omega \nabla \cdot (c_k\mb{q}_k ) v \diff x
    - \int_\Omega R_k v \diff x = 0 \quad \text{in}\quad\Omega.
\end{align}
The diffusion term can be treated exactly like in~\cite{arnold2002unified} with
the appropriate flux choices
\begin{align}
    \hat{c}_k = \avg{c_k} \text{ on } \Gamma_I,
    \quad \hat{c}_k = 0 \text{ on } \partial\Omega,
    \quad \hat{\sigma} = \avg{\nabla c_k} - D_k \delta_p \jump{c_k} \text{ on } \Gamma,
\end{align}
where $\delta_p$ is a penalty parameter that depends on the polynomial order $p$
of the finite elements. \rme{In \cite{hartmann2014higher}, it is given by}
\begin{equation}\label{eq:deltap}
    \delta_p = C_\mathrm{IP} \frac{p^2}{h},
\end{equation}
\rme{where $C_\mathrm{IP}>0$ needs to be large enough for stability, and $h$ is the cell diameter.}
The diffusion term results in
\begin{multline}
    -\int_\Omega \nabla \cdot (D_k \nabla c_k) v \diff x \approx \\
    \int_\Omega D_k \nabla c_k \cdot \nabla v \diff x
    - \int_\Gamma D_k (\jump{c_k} \cdot \avg{\nabla v} + \avg{\nabla c_k} \cdot \jump{v}) \diff s
    + \int_\Gamma D_k \delta_{p} \jump{c_k} \cdot \jump{v} \diff s.
\end{multline}
The advection-migration term combines the advection effect from the transport
velocity $\vec{v}$ and the migration effect from the potential gradient
$\nabla\Phi$, which can be treated with the combined variable $\vec{q}_k$. This
combined term is treated as in~\cite{cockburn2001runge} with the upwind flux, resulting in
\begin{multline}
    \int_\Omega \nabla \cdot (c_k\mb{q}_k ) v \diff x \approx \\
- \int_\Omega c_k \mb q_k \cdot \nabla v \diff x + \int_{\Gamma_I} \left(\avg{c_k \mb
q_k}+ \frac{|\mb q_k \cdot \mb n|}{2}\jump{c_k}\right)\cdot\jump{v}\diff s.
\end{multline}
Since the reaction term is cell-local, and it is assumed not to contain
derivatives of the primary solution variables, no further treatment is needed.

Note that the gradient potential term $\nabla \Phi$ in the mass conservation
equations \cref{eq:NP2} is only treated like a velocity in an advective term for
the concentrations. A penalty term is required for the stability of the DG
scheme and that is why we are using the CNP model. Indeed, the gradient
potential term in the charge conservation equation \cref{eq:charge2} will be
treated like an elliptic term, with appropriate penalization. Conversely, since
we have already penalized the concentration gradients $\nabla c_k$ in our
treatment of \cref{eq:NP2}, it is not required to do so for \cref{eq:charge2}.

We introduce the following short-hand notation for \cref{eq:charge2}
\begin{align}
    a_{km} &= F z_k(D_k-D_m),\\
    \kappa &= \sum_{k=1}^{m-1} z_k(z_k\mu_k - z_m\mu_m)F^2,\\
    f &= \sum_{k=1}^m z_k F R_k,
\end{align}
which results in
\begin{align}
    -\nabla \cdot \left[\sum_{k=1}^{m-1} a_{km} \nabla c_k\right] - \nabla
    \cdot\left( \kappa \nabla \Phi\right) - f = 0 \quad \text{in}\quad \Omega .
\end{align}
The fluxes for the concentration term are chosen as
\begin{align}
    \hat{c}_k = \avg{c_k} \text{ on } \Gamma_I,
    \quad \hat{c}_k = 0 \text{ on } \partial\Omega,
    \quad \hat{\bm\sigma} = \avg{\nabla c_k} \text{ on } \Gamma.
\end{align}
Without loss of generality, the case with only two species in the system ($k=2$) results in
\begin{align}
    \int_\Omega a_{1} \nabla c_1 \cdot \nabla v \diff x
    - \int_\Gamma a_1 (\jump{c_1} \cdot \avg{\nabla v} + \avg{\nabla c_1} \cdot \jump{v}) \diff s = \int_\Omega f\,v\diff x,
\end{align}
where $a_1 = a_{12}$. For the potential term, we choose fluxes with interior penalty
\begin{align}
    \hat{\Phi} = \avg{\Phi} \text{ on } \Gamma_I,
    \quad \hat{\Phi} = 0 \text{ on } \partial\Omega,
    \quad \hat{\bm\sigma} = \avg{\kappa \nabla \Phi} - \avg{\kappa} \delta_p \jump{\Phi} \text{ on } \Gamma.
\end{align}
This results in
\begin{multline}
    \int_\Omega \kappa \nabla \Phi \cdot \nabla v \diff x
    - \int_\Gamma (\jump{\Phi} \cdot \avg{\kappa \nabla v} + \avg{\kappa \nabla \Phi} \cdot \jump{v}) \diff s \\
    + \int_\Gamma \avg{\kappa} \delta_{p} \jump{\Phi} \cdot \jump{v} \diff s
    = \int_\Omega f\,v\diff x,
\end{multline}
with the standard symmetric interior
penalty flux. Note that in this example we assume the same test function space for all
concentrations $c_k$ and the potential $\Phi$. Therefore, we use the same
penalty parameter $\delta_p$. If the polynomial spaces have different orders,
the penalty parameters for both terms should be chosen accordingly.
\rme{In our numerical experiments, we use \cref{eq:deltap} with $C_\mathrm{IP}=10$}.

\subsection{Scheme}
\label{subsec:dg_scheme}
Let the short-hand notation for the current at the electrode boundary
\begin{align}
g &= \sum_{k=1}^{m} z_k F g_k,
\end{align}
where $g_k$ are the reaction rates from \cref{eq:electrode}.

The DG scheme for the problem
\cref{eq:NP2,eq:charge2,eq:inlet,eq:outlet,eq:electrode,eq:wall} is given by:
Find $c_k\in V_h^p$ for $k=1,\dots,m-1$, and $\Phi\in V_h^p$ such that
\begin{subequations}
\label{eq:dgscheme}
\begin{multline}
    \int_\Omega D_k \nabla c_k \cdot \nabla v_k \diff x
    - \int_{\Gamma_I} (\jump{c_k} \cdot \avg{D_k \nabla v_k} + \avg{D_k \nabla c_k} \cdot \jump{v_k}) \diff s
    \\
    + \int_{\Gamma_I} \avg{D_k} \delta_p \jump{c_k} \cdot \jump{v_k} \diff s
    - \int_\Omega c_k \mb q_k \cdot \nabla v_k \diff x \\
    + \int_{\Gamma_I} \left(\avg{c_k \mb
    q_k}+ \frac{|\mb q_k \cdot \mb n|}{2}\jump{c_k}\right)\cdot\jump{v_k}\diff s
    + \int_{\Gamma_\mathrm{out}} (\mb q_k \cdot \mb n)c_k v_k\diff s \\
    + \int_{\Gamma_\mathrm{in}} (\mb u \cdot \mb n ) c_k^\mathrm{in} v_k \diff s
    - \int_{\Gamma_\mathrm{elec}}  g_k v_k \diff s
    -\int_\Omega R_k v_k \diff x  = 0,\label{eq:NPDG}
\end{multline}
for all $v_k\in V_h^p$, for $k=1,\dots,m-1$, and
\begin{multline}
    \sum_{k=1}^{m-1}\left[\int_\Omega a_{km} \nabla c_k \cdot \nabla v_m \diff x
    - \int_{\Gamma_I} (\jump{c_k} \cdot \avg{a_{km} \nabla v_m} + \avg{a_{km} \nabla
    c_k} \cdot \jump{v_m}) \diff s \right] \\
    + \int_\Omega \kappa \nabla \Phi \cdot \nabla v_m \diff x
    - \int_{\Gamma_I} (\jump{\Phi} \cdot \avg{\kappa \nabla v_m}
    + \avg{\kappa \nabla \Phi} \cdot \jump{v_m}) \diff s \\
    + \int_{\Gamma_I} \avg{\kappa} \delta_p \jump{\Phi} \cdot \jump{v_m} \diff s
    - \int_{\Gamma_\mathrm{elec}} g v_m \diff s
    -\int_\Omega f v_m \diff x  = 0\numberthis, \label{eq:chargeDG}
\end{multline}
\end{subequations}
for all $v_m\in V_h^p$.

Recall that the eliminated concentration $c_m$ can be recovered from the
substitution \cref{eq:substitution}. Note that in \cref{eq:chargeDG}, the inlet
boundary conditions cancel due to electroneutrality.

\section{Numerical details and preconditioning}
\label{sec:pc}

A crucial aspect in the practical usability of this specific discretization is the
ability to solve it efficiently when applied to real world problems. Domains for
such applications are usually three-dimensional and nonlinear couplings appear
through species reactions. Consequently, the numerical solution approach has to
be scalable on a desired hardware architecture to provide feasible turnaround
times for experiments. As a first performance measure, we are focusing on
hexahedral mesh elements, which allow for a tensor product
element approach \cite{McRae2016}, a fact that we exploit in our implementation (note
that we are not exploiting Gauss-Legendre-Lobatto nodes, because we tend to use
higher-order quadrature rules compared to the polynomial degree).

After discretization, the nonlinear system arising
from~\cref{eq:dgscheme} has to be linearized and solved frequently. To solve the
nonlinear system, we employ a Newton method along with backtracking
linesearch~\cite{dennis1996numerical}. To ensure a scalable linear solving
procedure we have to employ appropriate preconditioning techniques \cite{wathen2015preconditioning}.

We can rewrite~\cref{eq:dgscheme} as a block residual equation
\begin{align*}
    \mathcal{R} =
    \begin{bmatrix}
        R_{\Phi}\\
        R_{c_1}\\
        \vdots \\
        R_{c_{m-1}}
    \end{bmatrix}
    = 0,
\end{align*}
where $R_\Phi$ is the residual of the discretized charge conservation equation \cref{eq:chargeDG}, and $R_{c_i}$, $i=1,\dots,m-1$, correspond to the residuals of the mass conservation equations \cref{eq:NPDG}.
Linearization of the residual $\mathcal{R}$ results in the block matrix
\begin{align*}
    J = \pder{\mathcal{R}}{(\Phi, c_1, \dots, c_{m-1})} &=
    \begin{bmatrix}
        \pder{R_{\Phi}}{\Phi} & \pder{R_{\Phi}}{c_1} & \dots & \pder{R_{\Phi}}{c_{m-1}} \\
        \pder{R_{c_1}}{\Phi} & \pder{R_{c_1}}{c_1} & \dots & \pder{R_{c_1}}{c_{m-1}} \\
        \vdots & \vdots & \ddots & \vdots \\
        \pder{R_{c_{m-1}}}{\Phi} & \pder{R_{c_{m-1}}}{c_1} & \dots & \pder{R_{c_{m-1}}}{c_{m-1}}
    \end{bmatrix}\\
    &=
    \begin{bmatrix}
        A_{\Phi \Phi} & A_{\Phi c_1} & \dots & A_{\Phi c_{m-1}}\\
        A_{c_1 \Phi} & A_{c_1 c_1} & \dots & A_{c_1 c_{m-1}}\\
        \vdots & \vdots & \ddots & \vdots \\
        A_{c_{m-1} \Phi} & A_{c_{m-1} c_1} & \dots & A_{c_{m-1} c_{m-1}}\\
    \end{bmatrix}.\numberthis \label{eq:blockmatrix}
\end{align*}
\rme{This block system is used in the following sections to determine preconditioning strategies.}

\subsection{Preconditioning for $p=1$}
\label{sec:pcp1}
\rme{We first describe the preconditioning strategy for a piecewise linear discretization.}

When dealing with block systems resulting from the discretizations of systems of PDEs with different underlying physics, it is
usually not sufficient and never optimal to use the same preconditioner for each
block. Analyzing the linearized block system \cref{eq:blockmatrix}, we immediately notice that
the~$A_{\Phi\Phi}$~block only consists of terms that resemble a Poisson equation,
for which preconditioners like Algebraic Multigrid (AMG) \cite{ruge1987algebraic}
combined with Conjugate Gradient (CG) iterations work well.

Continuing with the
diagonal, the~$A_{c_k c_k}$~blocks are not as simple and involve a diffusive term as well as
a (usually dominant) advective term with anisotropic coefficients.
In this case, \rme{classical} AMG will not work well since its heuristics require elliptic-like properties.
Instead, we apply preconditioned GMRES iterations \cite{saad1986gmres} to the $A_{c_k c_k}$ blocks,
considering two different preconditioners.
The first preconditioning approach is an additive Schwarz method (ASM) \cite{widlund1987additive}
where an incomplete LU factorization with 0 fill (ILU0) \cite{meijerink1977iterative} is applied to
each subdomain.
The subdomains are given by the parallel decomposition of the mesh.
The second approach is a Geometric Multigrid (GMG) \cite{brandt1977multi}
preconditioner where that same ASM-ILU0 method is used
as a smoother. The coarsest grid can be solved directly with a full LU factorization.
Although, naively solving the coarse grid this way does not scale to a large
number processors when parallel performance is desired. Because the coarse grid
and therefore the coarse discretization is usually small, even on computations
with hundreds of processors, the communication between the processors is the
bottleneck for performance. A solution is to redistribute and duplicate the
coarse grid matrix to predefined subsets of processors. This technique performs
redundant work but improves the overall solve time as shown in \cite{may2016extreme}.

We choose a block-triangular preconditioner of the form
\rme{
\begin{equation}
    \label{eq:P1}
    P =
    \begin{bmatrix}
        A_{\Phi \Phi} &  & & \\
        A_{c_1 \Phi} & A_{c_1 c_1} &  & \\
        \vdots & \vdots & \ddots &  \\
        A_{c_{m-1} \Phi} & A_{c_{m-1} A_{c_1}} & \dots & A_{c_{m-1} c_{m-1}}\\
    \end{bmatrix}^{-1}.
\end{equation}
}
Inverting a block-triangular matrix requires the inversion of the diagonal blocks,
which are approximated by
\begin{equation}
A^{-1}_{\Phi \Phi} \approx \mathrm{CG}(\mathrm{AMG}(A_{\Phi \Phi})),\qquad
A^{-1}_{c_k c_k} \approx \mathrm{GMRES}(\mathrm{ASM}(A_{c_k c_k})),
\end{equation}
in the first approach, and
\begin{equation}
A^{-1}_{\Phi \Phi} \approx \mathrm{CG}(\mathrm{AMG}(A_{\Phi \Phi})),\qquad
A^{-1}_{c_k c_k} \approx \mathrm{GMRES}(\mathrm{GMG}(A_{c_k c_k})),
\end{equation}
in the second.

\subsection{Preconditioning for high order}
\label{sec:pchigh}

\rme{AMG is more memory intensive for higher-order discretizations.
In this case, $p$-multigrid has been successful in the past \cite{fidkowski2005p,helenbrook2003analysis}.
For the $A_{\Phi\Phi}$ block, we consider a $p$-multigrid strategy where AMG is applied to a lower-order $p$-coarsened problem.
This coarsened problem uses the same DG discretization for $p=1$, where the penalization terms are appropriately updated for the lower order.
We include element-wise Jacobi as a relaxation scheme, where inter-element contributions are ignored to maintain locality.}

\rme{For the concentration blocks, we consider the same additive Schwarz method described for $p=1$. We do not consider the use of GMG for the high-order case.}

\rme{The block preconditioner \cref{eq:P1} with the described p-multigrid method for the potential block exhibits convergence issues for some numerical experiments of \cref{sec:results} for $p=3$. It is unclear to the authors why this is the case. Instead, a different ordering of the preconditioning steps works better: solving the concentration blocks first, then the potential blocks. We theorize that the couplings are important in the low-order coarse problem but have only heuristic evidence. We thus use the following preconditioner:}
\rme{
    \begin{equation}
    \label{eq:P2}
    P =
    \begin{bmatrix}
        A_{\Phi \Phi} & A_{\Phi c_1} & \dots & A_{\Phi c_{m-1}} \\
         & A_{c_1 c_1} &  & \\
         & \vdots & \ddots &  \\
         & A_{c_{m-1} A_{c_1}} & \dots & A_{c_{m-1} c_{m-1}}\\
    \end{bmatrix}^{-1}.
\end{equation}
}
\rme{Inverting this matrix only requires the inversion of the diagonal blocks,
which are approximated by
\begin{equation}
A^{-1}_{\Phi \Phi} \approx \mathrm{CG}(\mathrm{PMG}(A_{\Phi \Phi})),\qquad
A^{-1}_{c_k c_k} \approx \mathrm{GMRES}(\mathrm{ASM}(A_{c_k c_k})).
\end{equation}
}

\subsection{Implementation}

The implementation \rme{is part of the \texttt{EchemFEM} package \cite{echemfem}}. It is realized using \texttt{firedrake}~\cite{Rathgeber2016} as
a finite element discretization tool, which allows for automated code generation
and optimization. We leverage the capability to generate kernels of tensor
product finite elements~\cite{McRae2016} on quadrilateral~\cite{Homolya2016} and
extruded hexahedral meshes~\cite{Bercea2016}. Note that the algorithms and
discretization techniques presented here can also be applied to triangular and
tetrahedral elements. Furthermore, \texttt{firedrake} provides the necessary
infrastructure to realize the presented GMG preconditioner~\cite{Mitchell2016}.
The nonlinear and linear solves are handed off to the
\texttt{PETSc}~\cite{petsc-efficient} linear algebra backend, and \texttt{BoomerAMG}~\cite{henson2002boomeramg} from the
\texttt{hypre} library~\cite{hypre} is used for the AMG implementation.
\rme{For higher-order simulations, we use \texttt{firedrake}'s matrix-free capabilities \cite{kirby2018solver}.}
In~\cref{sec:appendix}
we provide all solver options, along with relative tolerances of each \rme{inner} Krylov
method, \rme{the outer FGMRES~\cite{saad1993flexible}}, and the outer Newton solver.

To help the convergence of Newton's method, we start with a constant initial guess for the concentrations using the inlet values $c_k^\mathrm{in}$. For the potential, we obtain an initial guess by first solving the charge conservation equation with constant concentrations (a Poisson equation). For this initial potential solution, we use the solver for the potential block defined in \cref{sec:appendix}, with a conjugate gradient tolerance of $10^{-2}$.

\section{Numerical results}
\label{sec:results}
In this section, we investigate the numerical performance of the DG scheme and preconditioning strategies.
We perform a convergence test for polynomial orders $p=1,2,3$ and strong scaling tests for a large parallel plate flow channel example for $p=1$ and $p=3$.
All tests are performed in three-dimensional domains.

\subsection{Convergence tests}
In \cref{sec:dg}, we established a discretization scheme based on well-known DG
methods. As a first step, we have to make sure that the presented method provides
the expected accuracy. Optimally, we would be able to prove an error estimate and
therefore a theoretical convergence rate through analysis of~\cref{eq:dgscheme}.

Another, although less thorough, way of determining the convergence rate is by
solving a manufactured problem and comparing the solution against the known
analytical solution. The way the scheme was set up, we expect the convergence
rate to be similar to the minimum rate of each individual
treated equation. Therefore, the convergence rate of the concentration variables
should be close to that of an advection-diffusion-reaction equation, i.e.
$p+\frac 1 2$ in the advection-dominated regime \cite{ayuso2009discontinuous},
and the convergence rate of the potential variable should be close to that of
a diffusion equation, i.e. $p+1$ \cite{arnold2002unified}.

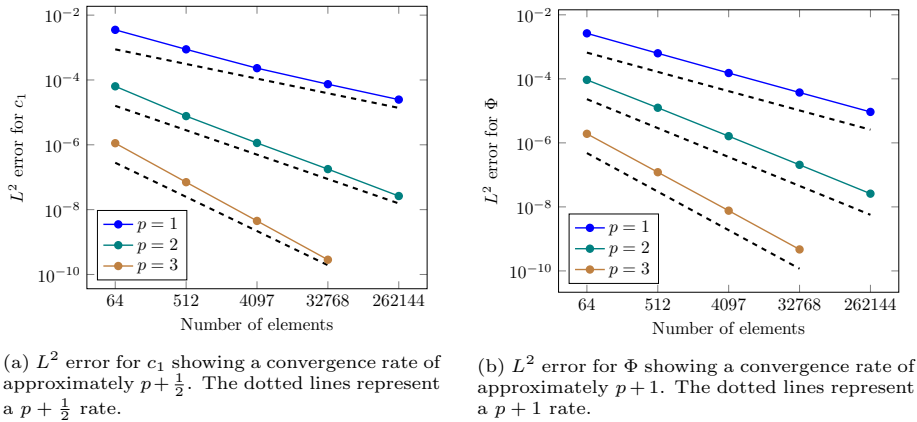
\begin{figure}[htb!]
\centering
\begin{subfigure}[b]{0.48\linewidth}
\begin{tikzpicture}[scale=0.66]
\begin{axis}[xmode=log,ymode=log,xminorticks=false,yminorticks=false,
    xlabel={Number of elements},
    ylabel={$L^2$ error for $c_1$},
    xtick={64,512,4096,32768,262144},
    xticklabels={64,512,4097,32768,262144},
    legend pos = south west,
    ]
    \addplot[thick,color=blue,mark=*]
    table [x=ncells,y=l2error,col sep=comma] {convergence_p1_c.csv};
    \addplot[thick,color=teal,mark=*]
    table [x=ncells,y=l2error,col sep=comma] {convergence_p2_c.csv};
    \addplot[thick,,color=brown,mark=*]
    table [x=ncells,y=l2error,col sep=comma] {convergence_p3_c.csv};
    \addplot[very thick,color=black,dashed] coordinates {
        (64,8.800E-04)
        (512,3.111E-04)
        (4096,1.100E-04)
        (32768,3.889E-05)
        (262144,1.375E-05)
    };
    \addplot[very thick,color=black,dashed] coordinates {
        (64,1.593E-05)
        (512,2.815E-06)
        (4096,4.977E-07)
        (32768,8.797E-08)
        (262144,1.555E-08)
    };
    \addplot[very thick,color=black,dashed] coordinates {
        (64,2.800E-07)
        (512,2.475E-08)
        (4096,2.187E-09)
        (32768,1.933E-10)
    };
    \legend{$p=1$,$p=2$,$p=3$};
\end{axis}
\end{tikzpicture}
\caption{$L^2$ error for $c_1$ showing a convergence rate of approximately $p+\frac 12$. The dotted lines represent a $p+\frac 12$ rate. }
\end{subfigure}
\hfill
\begin{subfigure}[b]{0.48\linewidth}
\begin{tikzpicture}[scale=0.66]
\begin{axis}[xmode=log,ymode=log,xminorticks=false,yminorticks=false,
    xlabel={Number of elements},
    ylabel={$L^2$ error for $\Phi$},
    xtick={64,512,4096,32768,262144},
    xticklabels={64,512,4097,32768,262144},
    legend pos = south west,
    ]
    \addplot[thick,color=blue,mark=*]
    table [x=ncells,y=l2error,col sep=comma] {convergence_p1_u.csv};
    \addplot[thick,color=teal,mark=*]
    table [x=ncells,y=l2error,col sep=comma] {convergence_p2_u.csv};
    \addplot[thick,,color=brown,mark=*]
    table [x=ncells,y=l2error,col sep=comma] {convergence_p3_u.csv};
    \addplot[very thick,color=black,dashed] coordinates {
        (64,6.600E-04)
        (512,1.650E-04)
        (4096,4.125E-05)
        (32768,1.031E-05)
        (262144,2.578E-06)
    };
    \addplot[very thick,color=black,dashed] coordinates {
        (64,2.318E-05)
        (512,2.897E-06)
        (4096,3.621E-07)
        (32768,4.526E-08)
        (262144,5.658E-09)
    };
    \addplot[very thick,color=black,dashed] coordinates {
        (64,4.800E-07)
        (512,3.000E-08)
        (4096,1.875E-09)
        (32768,1.172E-10)
    };
    \legend{$p=1$,$p=2$,$p=3$};
\end{axis}
\end{tikzpicture}
\caption{$L^2$ error for $\Phi$ showing a convergence rate of approximately $p+1$. The dotted lines represent a $p+1$ rate. }
\end{subfigure}
\caption{Convergence tests for the manufactured solution.}
\label{fig:convergence_rate}
\end{figure}

We use the nondimensional version of a two-ion system
with inverse P\'eclet numbers $\hat{D}_1=5\times 10^{-6}$ and $\hat{D}_2=10^{-5}$, and charge
numbers $z_1=2$ and $z_2=-2$.
The analytical manufactured solution is chosen as
\begin{align}
    c_{1\text{ex}} &= \cos{x} + \sin{y} + 3,\\
    \Phi_{\text{ex}} &= \sin{x} + \cos{y} + 3,
\end{align}
on the domain $\Omega \in (0,1)^3$. All boundary conditions on $\partial\Omega$
are prescribed using the exact solutions. Additionally, a
parabolic velocity profile is defined with $\vec{u} = (6 y(1-y), 0, 0)$ in the whole domain.
We choose an initial mesh size using 64 hexahedral elements with a polynomial
order of $p=1$ and compute the $L^2$ error norms of each variable
\begin{align}
   \|c_1 - c_{1\text{ex}}\|_2, \qquad  \|\Phi - \Phi_{\text{ex}}\|_2.
\end{align}
Then, we subsequently halve the element size $h_e$ using uniform refinement. The
process is repeated for polynomial orders $p=2$ and $p=3$.
\Cref{fig:convergence_rate} shows the results for this test and one can observe
a convergence rate of approximately $p+1$ for $\Phi$. For $c_1$, the rate appears
to be $p+1$ for the coarser meshes, until the later refinements reveals the asymptotic
rate of $p+\frac 12$. This asymptotic rate is not reached for $p=3$.

Looking at these results, one has to keep in mind that we determined the
convergence rates on a Cartesian grid, and they might deteriorate on unstructured
or curved grids~\cite{cockburn2008optimal,richter1988optimal}.

\rme{
We now briefly investigate how the accuracy of the method relates to computational cost.
For the manufactured solution problem described above, we calculate the total computational time and total $L^2$ error as we vary the size of the problem.
The total error is the sum of the concentration and potential errors.
On an Intel 18-core Xeon E5-2695 v4 architecture with 2.1 GHz clock speed, we use two nodes and a total of 64 processors.
For $p=1$, we use the preconditioner described in \cref{sec:pcp1} with the additive Schwarz method for the concentration blocks, and for $p=3$ we use the preconditioner for high order described in \cref{sec:pchigh}.
}

\begin{table}[htb!]
    \centering
    \caption{Computational time and total $L^2$ errors for varying problem sizes of the toy problem.}
    \label{tab:cpu}
    \begin{tabular}{lcccc}
        \toprule
        &  \multicolumn{2}{c}{$p=1$} & \multicolumn{2}{c}{$p=3$} \\
        \midrule
        \# DoFs & CPU time (s) & $L^2$ error & CPU time (s) & $L^2$ error \\
        \midrule
        65,536	& 15.31 & \num{3.84E-04} & 29.92 &	\num{1.92E-07} \\
        524,288	& 18.17 & \num{1.14E-04} & 41.18 & \num{1.12E-08} \\
        4,194,304 &	34.95 &	\num{3.40E-05} & 111.3 & \num{7.61E-10} \\
        33,554,432 & 209.1 & \num{9.82E-06} & - & - \\
        \bottomrule
    \end{tabular}
\end{table}

\rme{
In \cref{tab:cpu}, we provide computational timings and total $L^2$ errors for $p=1$ and $p=3$ as we vary the problem size.
We note that the computational resources are fixed so there is insufficient memory for the largest high-order case.
Additionally, the first row of results is around a thousand degrees of freedom per processor, which is not typically cost-efficient.
Even so, the smallest $p=3$ problem is twice as accurate and seven times faster than the largest $p=1$ problem.
It is evident that increasing polynomial order is preferable to mesh refinement when seeking high accuracy results.
}

\subsection{Parallel plate flow channel}
\label{sec:scaling}
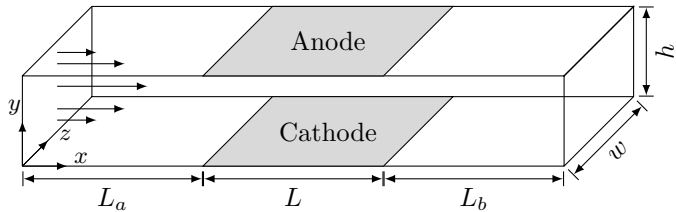
\begin{figure}[htb!]
    \centering
        \begin{tikzpicture}[>=latex,scale=1.2]
            \pgfmathsetmacro{\La}{2}
            \pgfmathsetmacro{\L}{2}
            \pgfmathsetmacro{\h}{1}
            \pgfmathsetmacro{\hz}{2}
            \pgfmathsetmacro{\x}{2*\La+\L}
            \pgfmathsetmacro{\y}{\h}
            \pgfmathsetmacro{\z}{\hz}
            \path (0,0,\z) coordinate (A) (\x,0,\z) coordinate (B) (\x,0,0) coordinate (C) (0,0,0)
            coordinate (D) (0,\y,\z) coordinate (E) (\x,\y,\z) coordinate (F) (\x,\y,0) coordinate (G)
            (0,\y,0) coordinate (H);
            \draw (A)--(B)--(C)--(G)--(F)--(B) (A)--(E)--(F)--(G)--(H)--(E);
            \draw (A)--(D)--(C) (D)--(H);
            \coordinate (A1) at (\La,0,\z);
            \coordinate (A2) at (\La+\L,0,\z);
            \coordinate (A3) at (\La+\L,0,0);
            \coordinate (A4) at (\La,0,0);
            \coordinate (B1) at (\La,\y,\z);
            \coordinate (B2) at (\La+\L,\y,\z);
            \coordinate (B3) at (\La+\L,\y,0);
            \coordinate (B4) at (\La,\y,0);
            \coordinate (V1) at (0,\y/2,\z/2);
            \coordinate (V2) at (\La/2,\y/2,\z/2);
            \coordinate (W1) at (0,\y/4,\z/2);
            \coordinate (W2) at (\La*3/8,\y/4,\z/2);
            \coordinate (U1) at (0,3*\y/4,\z/2);
            \coordinate (U2) at (\La*3/8,3*\y/4,\z/2);
            \coordinate (W1a) at (0,\y/8,\z/2);
            \coordinate (W2a) at (\La*7/32,\y/8,\z/2);
            \coordinate (U1a) at (0,7*\y/8,\z/2);
            \coordinate (U2a) at (\La*7/32,7*\y/8,\z/2);
            \draw[->] (V1) -- (V2);
            \draw[->] (W1) -- (W2);
            \draw[->] (U1) -- (U2);
            \draw[->] (W1a) -- (W2a);
            \draw[->] (U1a) -- (U2a);

            \draw[fill=gray!30] (A1)--(A2)--(A3)--(A4)--cycle;
            \draw[fill=gray!30] (B1)--(B2)--(B3)--(B4)--cycle;
            \node[] at ($(A1)!0.5!(A3)$) {Cathode};
            \node[] at ($(B1)!0.5!(B3)$) {Anode};
            \draw[thin,|<->|] ($(A)+(0,-4pt)$) -- node[below]{$L_a$}($(A1)+(0,-4pt)$);
            \draw[thin,|<->|] ($(A1)+(0,-4pt)$) -- node[below]{$L$}($(A2)+(0,-4pt)$);
            \draw[thin,|<->|] ($(A2)+(0,-4pt)$) -- node[below]{$L_b$}($(B)+(0,-4pt)$);
            \draw[thin,|<->|] ($(C)+(4pt,0)$) -- node[below,sloped]{$h$}($(G)+(4pt,0)$);
            \draw[thin,|<->|] ($(B)+(-45:4pt)$) -- node[below,sloped]{$w$}($(C)+(-45:4pt)$);
            \coordinate (O) at (0,0,2);
            \draw[->] (O) --+ (0.5,0,0) node[anchor=west,yshift=3pt,xshift=-1pt]{\small $x$};
            \draw[->] (O) --+ (0,0.5,0) node[anchor=south, xshift=-3pt,yshift=-1pt]{\small $y$};
            \draw[->] (O) --+ (0,0,-0.75) node[anchor=south west, xshift=0.5pt,yshift=-3.5pt]{\small $z$};
        \end{tikzpicture}
        \caption{Parallel plate reactor diagram.}\label{fig:reactor}
\end{figure}

We consider a 3D extension of the parallel plate reactor from
\cite{bortels1996multi}. A schematic of the reactor is given in
\cref{fig:reactor}. The distance between the two electrodes is $h=\SI{0.01}{\metre}$,
and their length is $L=\SI{0.02}{\metre}$. The inlet and outlet regions included in the
domain are of length $L_a=L_b = \SI{0.05}{\metre}$, and the channel has a width $w=\SI{0.06}{\metre}$.

The velocity field $\mb u = (u_x,u_y,u_z)$ in the reactor is given by the parabolic profile
\begin{equation}
    u_x = \frac{6 u_\mathrm{avg}}{h^2} y (h -y),\quad u_y = u_z = 0.
\end{equation}
Choosing zero-Neumann boundary conditions on $xy$-faces means that the channel is periodic in the $z$-direction.
For this test case, we choose $u_\mathrm{avg} = \SI{0.03}{\metre\per\second}$.

We consider the three-ion system from \cite{bortels1996multi} with an electrolyte solution consisting of 0.01 M \chem{CuSO_{4}} and 1.0 M \chem{H_2SO_4}. The ions are \chem{H^+}, \chem{Cu^{2+}}, and \chem{SO_4^{2-}} with properties given in \cref{tab:ionparams}.
\rone{This results in average P\'eclet numbers, $u_\mathrm{avg} \frac{h}{D_{k}}$, of $3.22\times 10^4$, $4.17\times 10^5$, and $2.82\times 10^5$ for \chem{H^+}, \chem{Cu^{2+}}, and \chem{SO_4^{2-}}, respectively.}

\begin{table}[htb!]
    \centering
    \caption{Physical parameters of the ions for the parallel plate flow reactor.}\label{tab:ionparams}
    \begin{tabular}{cccc}
        \toprule
        Ion              & Diffusion coefficient (\si{m^2 s^{-1}}) & Inlet concentration (\si{M}) & charge \\
        \midrule
        \chem{Cu^{2+}}   & $7.20 \times 10^{-10}$                   & 0.01                         & 2      \\
        \chem{SO_4^{2-}} & $10.65 \times 10^{-10}$                  & 1.01                         & -2     \\
        \chem{H^+}       & $93.12 \times 10^{-10}$                  & 2.0                          & 1      \\
        \bottomrule
    \end{tabular}
\end{table}

At the electrodes, the electrochemical reaction is given by the Butler-Volmer
equation \cref{eq:bv} for the reacting ion, \chem{Cu^{2+}}, and no reaction for
the other ions. The parameters for the Butler-Volmer equation are $\gamma=1$, $\alpha_1 = \alpha_2 = 0.5$, $T=\SI{298.15}{\kelvin}$,
$n=2$, and $c_{\chem{Cu^{2+}}}^*$ is given by the inlet concentration of \chem{Cu^{2+}} in
\cref{tab:ionparams}. The applied potential $\Phi_\mathrm{app}$ is zero at the
anode and \SI{0.03}{\volt} at the cathode.
For this redox reaction, \chem{Cu^{2+}} is the oxidant, and the solid copper of the electrode is the reductant (with constant concentration).
To introduce variations in the $z$ direction, we consider the following parabolic profile for the exchange current density:
\begin{equation}
    J_0 = \frac{3}{5} \bar J_0 \left[2 - \left(\frac{z - w/2}{w/2}\right)^2\right],
\end{equation}
which has an average of $\bar J_0=\SI{30}{\ampere\per\square\metre}$ over the width of the channel.

\begin{figure}[htb!]
    \centering
    \begin{tikzpicture}[>=latex]
    \tikzstyle{every node}=[font=\small]
    \node[inner sep=0pt] (mesh) at (0,0)
        {\includegraphics[width=310pt]{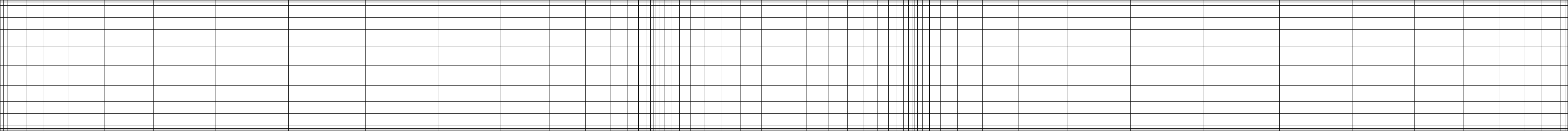}};
    \node[inner sep=0pt] (zoom) at (0,-2.4)
        {\includegraphics[width=140pt]{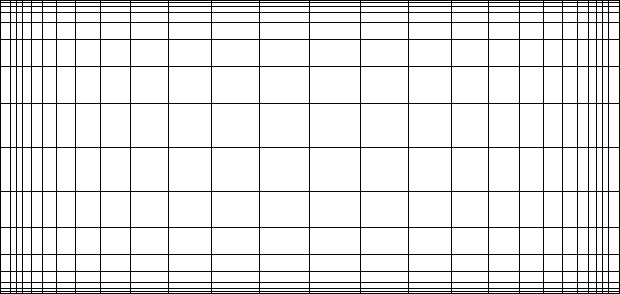}};
    \draw (-.97,-.5) -- (0.97,-.5) -- (0.97,0.5) -- (-.97,0.5) -- cycle;
    \draw[->,thick] (-.97,-.5) -- (zoom.north west);
    \draw[->,thick] (.97,-.5) -- (zoom.north east);
    \node at (0,.7) {Anode};
    \node at (0,-.7) {Cathode};
    \end{tikzpicture}
    \caption{Cross-section of the coarse mesh in the $xy$-plane and zoom-in of the electrode region.}\label{fig:mesh}
\end{figure}
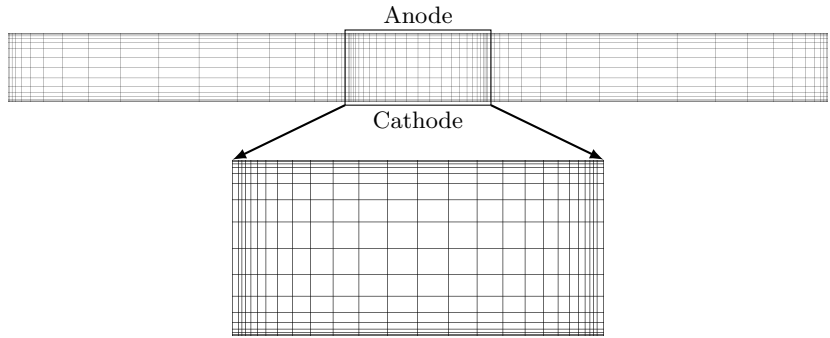

We \rme{first} use a structured hexahedral mesh such that, in the $xy$-plane, elements are concentrated in the electrode region as illustrated in \cref{fig:mesh}, and elements are equally spaced in the $z$-plane.
The coarse mesh of $64\times 16 \times 8$ elements is uniformly refined three times, resulting in a mesh with \rme{4,194,304} elements.
\rme{Note that this mesh has elements with large aspect ratios.}

\begin{figure}
    \centering
    \includegraphics[width=310pt]{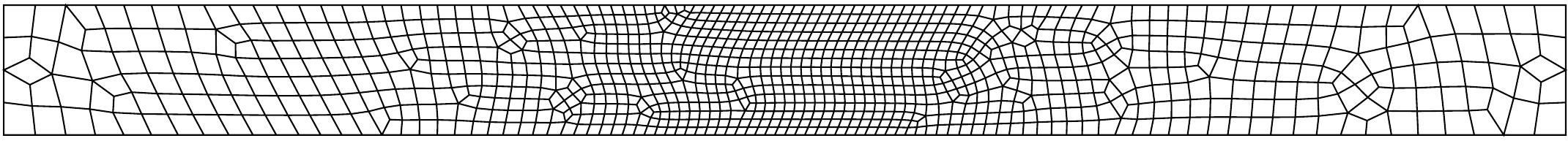}
    \caption{Cross-section of the unstructured coarse mesh in the $xy$-plane.}\label{fig:mesh2}
\end{figure}

\rtwo{
We also consider a mesh where the $xy$-plane is an unstructured mesh as illustrated in \cref{fig:mesh2}.
Similarly to the structured mesh, the coarse 2D mesh is extruded in the $z$-direction with eight equally-spaced elements.
After uniformly refining three times, we obtain a mesh with 4,153,344 elements.
}

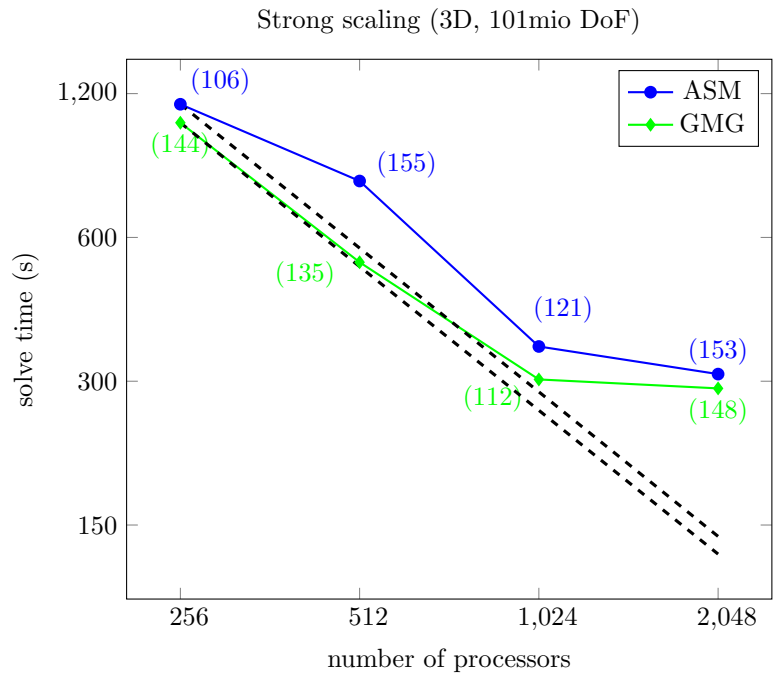
\begin{figure}[htb!]
\centering
\begin{tikzpicture}
    \begin{loglogaxis}[
      log ticks with fixed point,
      scaled ticks=false,
      tick scale binop=\times,
        width=4in,
        xlabel={number of processors},
        ylabel={solve time (s)},
        title={Strong scaling (3D, 101mio DoF)},
        xtick={256,512,1024,2048},
        ytick={1200,600,300,150},
        log basis x=2,
        xticklabel={
            \pgfkeys{/pgf/fpu=true}
            \pgfmathparse{int(2^int(\tick))}
            \pgfmathprintnumber[fixed]{\pgfmathresult}
        },
        ]

        \addplot[thick,color=blue,mark=*] coordinates {
          (256,1140)
          (512,787.6)
          (1024,354.7)
          (2048,310.5)
            }
            node[pos=0.0, above right]{(106)}
            node[pos=0.333, above right]{(155)}
            node[pos=0.666, above right]{(121)}
            node[pos=1.0, above]{(153)}
            ;

        \addplot[thick,color=green,mark=diamond*] coordinates {
            (256,1042.8)
            (512,532.3)
            (1024,302.6)
            (2048,289.8)
            }
            node[pos=0.0, below]{(144)}
            node[pos=0.33, below left]{(135)}
            node[pos=0.69, below left]{(112)}
            node[pos=1.0, below]{(148)}
                ;
        \addplot[very thick,color=black,dashed] coordinates {
          (256,1140)
          (512,570)
          (1024,285)
          (2048,142)
        };
        \addplot[very thick,color=black,dashed] coordinates {
          (256,1042.8)
          (512,521.4)
          (1024,260.7)
          (2048,130.4)
        };
        \legend{ASM,GMG};
    \end{loglogaxis}
    \end{tikzpicture}
    \caption{Strong scaling timings for the parallel flow plane channel for $p=1$. The total number of outer GMRES iterations are given in parentheses.}
    \label{fig:scaling}
\end{figure}

We now investigate the strong scaling of the preconditioners for the parallel flow plane channel.
The experiments are performed on an Intel 18-core Xeon E5-2695 v4 architecture with 2.1 GHz clock speed.
Using elements of $p=1$ polynomial order results in 100,663,296 degrees of freedom
\rme{for the structured mesh, and 99,680,256 degrees of freedom for the unstructured mesh.}

\begin{figure}[htb!]
    \centering
    \begin{tikzpicture}
        \begin{loglogaxis}[
          log ticks with fixed point,
          scaled ticks=false,
          tick scale binop=\times,
            width=4in,
            xlabel={number of processors},
            ylabel={solve time (s)},
            title={Strong scaling (3D unstructured, 100mio DoF)},
            xtick={256,512,1024,2048},
            ytick={640,320,160,80},
            log basis x=2,
            xticklabel={
                \pgfkeys{/pgf/fpu=true}
                \pgfmathparse{int(2^int(\tick))}
                \pgfmathprintnumber[fixed]{\pgfmathresult}
            },
            ]
    
            \addplot[thick,color=blue,mark=*] coordinates {
              (256,479.1)
              (512,282.1)
              (1024,162.8)
              (2048,126.5)
                }
                node[pos=0.0, below left]{(82)}
                node[pos=0.333, below left]{(86)}
                node[pos=0.666, below]{(74)}
                node[pos=1.0, below]{(80)}
                ;
    
            \addplot[thick,color=green,mark=diamond*] coordinates {
                (256,631.8)
                (512,346.4)
                (1024,207.5)
                (2048,187.1)
                }
                node[pos=0.0, above right]{(80)}
                node[pos=0.33, above right]{(75)}
                node[pos=0.69, above right]{(72)}
                node[pos=1.0, above]{(80)}
                    ;
            \addplot[very thick,color=black,dashed] coordinates {
              (256,479.1)
              (512,239.6)
              (1024,119.8)
              (2048,59.9)
            };
            \addplot[very thick,color=black,dashed] coordinates {
              (256,631.8)
              (512,315.9)
              (1024,158)
              (2048,79)
            };
            \legend{ASM,GMG};
        \end{loglogaxis}
        \end{tikzpicture}
        \caption{Strong scaling timings for the parallel flow plane channel for $p=1$. The total number of outer GMRES iterations are given in parentheses.}
        \label{fig:scaling2}
\end{figure}
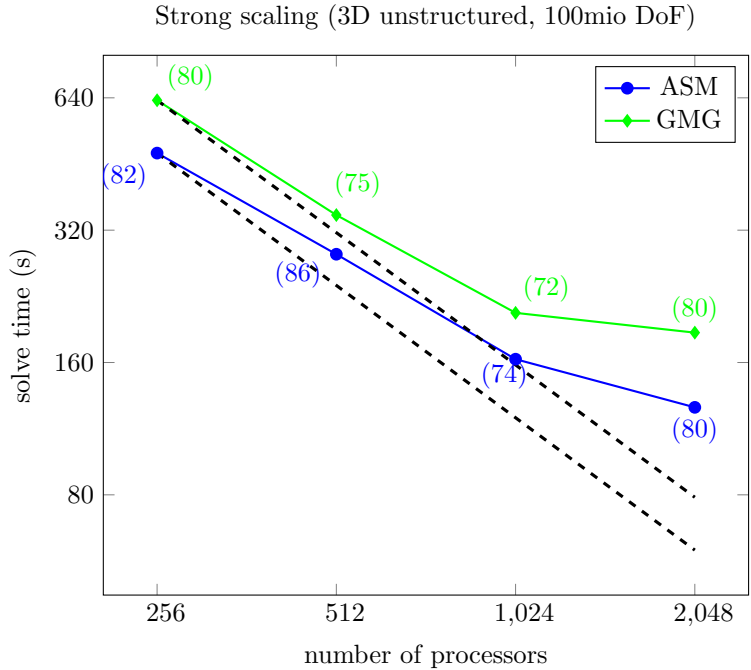

We begin with 256 processors (on 8 nodes) and successively double the number of processors and nodes.
In \cref{fig:scaling}, we illustrate the strong scaling results for both preconditioning approaches described in \cref{sec:pc}\rme{on the structured mesh described above, and in \cref{fig:scaling2} on the unstructured mesh}.
ASM refers to using the additive Schwarz method with ILU0 on the subdomains as a preconditioner for the concentration blocks, while GMG refers to using a geometric multigrid method with that same ASM as a smoother.
For both ASM and GMG, we present solution times and the total number of outer \rme{FGMRES} iterations.
In each case, four Newton iterations are required to reach the prescribed tolerance.

\rme{First, we discuss the scaling results for the unstructured mesh in \cref{fig:scaling}.}
We observe that, overall, ASM is the more expensive method.
Further investigation of the iteration counts of the inner GMRES used for the concentration blocks reveals that GMG only takes a few iterations, while ASM often takes around 50-100 iterations (sub-iterations not shown here).
Now, in terms of scaling, ASM scales almost linearly until 2,048 processors.
However, it is expected that ASM becomes a weaker preconditioner as the number of subdomains increases.
For GMG, the scaling is close to being linear form 256 to 1,024 processors.
Naturally, the direct solver used for the coarse problem of GMG does not scale linearly and eventually dominates the computational cost at 2,048 processors.
Note that there is a decrease in iterations at 512 and 1,024 processors, which may be due to how the mesh is partitioned, thus increasing the performance of the smoother and the apparent quality of the scaling.
Recall that we do not have control over the domain decomposition, which in general will not be ideal.
In particular, the criticality of the electrode regions should be considered in a more controlled mesh partition.

\rtwo{
Now we compare the scaling results for the unstructured mesh in \cref{fig:scaling2}.
We observe that both methods have similar almost-linear scaling, with a bump in cost at 2,048 processors.
Differently from the structured mesh case, ASM is the least expensive method.
In fact, the GMRES iteration count for the concentration blocks for ASM are lower in this case, typically less than 30 (sub-iterations not shown here).
The structured mesh has some elements with very large aspect ratio.
This could affect the conditioning of the mesh, resulting in larger number of iterations.
We note that the outer GMRES iterations is also lower for the unstructured case.
}

\rme{
Finally, we look at the scaling results for $p=3$ and the preconditioning strategy described in \cref{sec:pchigh}.
We note that the $p$-multigrid method used for the potential block performs poorly on meshes with high element aspect ratios such as the ones used above.
The scaling tests are therefore performed on a mesh with moderate aspect ratios.
The width of the channel is changed to $w=\SI{0.01}{\metre}$.
We consider an unstructured coarse mesh similar to the one illustrated in \cref{fig:mesh2} but with slightly larger elements and extrude it in the $z$-direction with 6 elements.
For $p=3$, we uniformly refine the mesh 2 times for a total of 254,208 elements.
One more refinement gives 2,033,664 elements, which is used for $p=1$.
Both problems consist of 48,807,936 degrees of freedom.
}

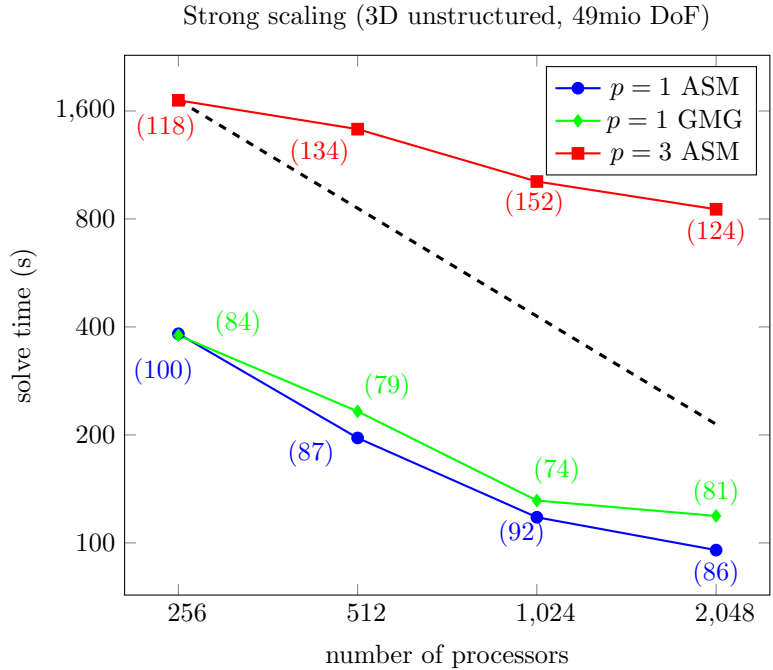
\begin{figure}[htb!]
    \centering
    \begin{tikzpicture}
        \begin{loglogaxis}[
          log ticks with fixed point,
          scaled ticks=false,
          tick scale binop=\times,
            width=4in,
            xlabel={number of processors},
            ylabel={solve time (s)},
            title={Strong scaling (3D unstructured, 49mio DoF)},
            xtick={256,512,1024,2048},
            ytick={1600,800,400,200,100},
            log basis x=2,
            xticklabel={
                \pgfkeys{/pgf/fpu=true}
                \pgfmathparse{int(2^int(\tick))}
                \pgfmathprintnumber[fixed]{\pgfmathresult}
            },
            ]
            \addplot[thick,color=blue,mark=*] coordinates {
                (256,383)
                (512,196.2)
                (1024,117.9)
                (2048,95.5)
                  }
                  node[pos=0.05, below left]{(100)}
                  node[pos=0.333, below left]{(87)}
                  node[pos=0.666, below]{(92)}
                  node[pos=1.0, below]{(86)}
                  ;
            \addplot[thick,color=green,mark=diamond*] coordinates {
                (256,379)
                (512,232.7)
                (1024,131.3)
                (2048,118.9)
                  }
                  node[pos=0.05, above right]{(84)}
                  node[pos=0.333, above right]{(79)}
                  node[pos=0.666, above right]{(74)}
                  node[pos=1.0, above]{(81)}
                  ;
            \addplot[thick,color=red,mark=square*] coordinates {
              (256,1713.9)
              (512,1424.0)
              (1024,1017.5)
              (2048,851.4)
                }
                node[pos=0.05, below left]{(118)}
                node[pos=0.333, below left]{(134)}
                node[pos=0.666, below]{(152)}
                node[pos=1.0, below]{(124)}
                ;
            \addplot[very thick,color=black,dashed] coordinates {
              (256,1713.9)
              (512,857)
              (1024,428.5)
              (2048,214.2)
            };
            \legend{$p=1$ ASM, $p=1$ GMG, $p=3$ ASM};
        \end{loglogaxis}
        \end{tikzpicture}
        \caption{Strong scaling timings for the parallel flow plane channel with ASM. The total number of outer GMRES iterations are given in parentheses.}
        \label{fig:scaling3}
\end{figure}
\rme{
We display the high-order scaling results in \cref{fig:scaling3}.
We observe that the scaling for $p=3$ is about order 0.5, i.e. the cost is halved after quadrupling the number of processors.
The computational timings for the smallest case is more than four times higher for $p=3$ than $p=1$.
Nevertheless, this figure does not take into account the numerical accuracy of the results.
Indeed, we recall from \cref{tab:cpu} that higher accuracy can be obtained for higher-order problems of much smaller sizes.
Strong scaling for $p=1$ as at best linear, but increasing the polynomial order results in an exponential increase in accuracy.
}

\section{Conclusion}

In this work, we propose a high-order Discontinuous Galerkin (DG) scheme and
block preconditioners for the electroneutral Nernst-Planck equations. We choose
a model formulation where the electroneutrality condition is substituted to
obtain a charge conservation equation. In the DG scheme, this equation is
treated as an elliptic equation for the electrical potential, whereas in the
mass conservation equations, an upwind scheme is used for the combined
advection-migration term. In numerical convergence tests, the scheme exhibits
a $p+1$ order of convergence for the potential and a $p+\frac 12$ order of convergence
for the concentrations.

\rme{For $p=1$,}
both proposed block preconditioners use Algebraic Multigrid for the potential block whereas an Additive Schwarz method with ILU is used for the concentration blocks, by itself or as a smoother for a Geometric Multigrid method (GMG). While both methods exhibit decent strong scaling, the GMG-based approach is less expensive \rme{on the structured mesh}, although its scaling is limited by the direct solver used for the coarse grid solution.
\rme{
For higher order, a $p$-multigrid method is used for the potential block, where Algebraic Multigrid is applied to the lower-order coarsened problem.
The strong scaling is less than linear in the high-order case.
 }

For future extensions, coupling this solver with fluid dynamics is of particular interest.
The Stokes or incompressible Navier-Stokes equations are relevant for flow reactors as well as the pore-scale of porous electrode systems.
On the larger length scale, Darcy flow is appropriate for porous electrodes, and its multiphase extension for gas-diffusion electrodes (GDE).

Further work is needed for preconditioning in the case where electrolyte bulk reactions are present, where reaction rate constants can vary wildly.
Porous electrodes are also interesting since cathodic/anodic reactions occur throughout the domain rather than as boundary conditions.

\section*{Acknowledgments}
The authors would like to thank Pablo Brubeck for his help on implementing $p$-multigrid.
This work was performed under the auspices of the U.S. Department of Energy by Lawrence Livermore National Laboratory under Contract DE-AC52-07NA27344 and was supported by the LLNL-LDRD program under project numbers 19-ERD-035 and 22-SI-006. LLNL Release Number LLNL-JRNL-829967.

\bibliography{bibfile}

\appendix

\def\appendixname{Appendix}
\section{Solver parameters}
\label{sec:appendix}

Here we provide PETSc options used to construct the solvers described in \cref{sec:pc}.
The parameters used for the $p=1$ case are provided in \cref{lst:petsc_options}, and for $p>1$ in \cref{lst:petsc_options_higherorder}.

\begin{lstlisting}[
    label={lst:petsc_options},
    language=Python,
    frame=singleo ,
    basicstyle=\footnotesize\ttfamily,
    numbers=left,
    numbersep=5pt,
    caption={Solver options for the discretized electrochemical system for $p=1$.}
]
"snes_type": "newton",
"snes_rtol": 1E-6,
"ksp_type": "fgmres",
"ksp_rtol": 1E-3,
"pc_type": "fieldsplit",
% For the potential block
"fieldsplit_0": {
    "ksp_rtol": 1E-1,
    "ksp_type": "cg",
    "pc_type": "hypre",
    "pc_hypre_boomeramg": {
            "strong_threshold": 0.7,
            "coarsen_type": "HMIS",
            "agg_nl": 3,
            "interp_type": "ext+i",
            "agg_num_paths": 5,
    },
},
% For each concentration block k, using the ASM preconditioner:
"fieldsplit_[k]": {
    "ksp_rtol": 1E-1,
    "ksp_type": "gmres",
    "pc_type": "asm",
    "sub_pc_type": "ilu",
    },
% Otherwise, for each concentration block k, using the GMG preconditioner:
"fieldsplit_[k]": {
    "ksp_rtol": 1E-1,
    "ksp_type": "gmres",
    "pc_type": "mg",
    "mg_levels_ksp_type": "richardson",
    "mg_levels_pc_type": "asm",
    "mg_levels_sub_pc_type": "ilu",
    "mg_coarse": {
        "pc_type": "python",
        "pc_python_type": "firedrake.AssembledPC",
        "assembled": {
                "mat_type": "aij",
                "pc_type": "telescope",
                "pc_telescope_reduction_factor": REDFACTOR,
                "pc_telescope_subcomm_type": "contiguous",
                "telescope_pc_type": "lu",
                "telescope_pc_factor_mat_solver_type": "mumps"
        }
    },
}
\end{lstlisting}

\begin{lstlisting}[
    label={lst:petsc_options_higherorder},
    language=Python,
    frame=singleo ,
    basicstyle=\footnotesize\ttfamily,
    numbers=left,
    numbersep=5pt,
    caption={Solver options for the discretized electrochemical system for $p>1$.}
]
"mat_type": "matfree",
"snes_type": "newton",
"snes_rtol": 1E-6,
"ksp_type": "fgmres",
"ksp_rtol": 1E-3,
"pc_type": "fieldsplit",
% For the potential block
"fieldsplit_[m]": {
    "ksp_rtol": 1E-1,
    "ksp_type": "cg",
    "pc_type": "python",
    "mat_type": "matfree",
    "pc_python_type": __name__+".CoarsenPenaltyPMGPC",
    "pmg_mg_levels_ksp_type": "chebyshev",
    "pmg_mg_levels_ksp_max_it": 4,
    "pmg_mg_levels_ksp_norm_type": "unpreconditioned",
    "pmg_mg_levels_": {
        "pc_type": "python",
        "pc_python_type": __name__ + "." + "CellIntegralPC",
        "assembled_pc_type": "jacobi",
    },
    "pmg_coarse_mat_type": "aij",
    "pmg_mg_coarse": {
        "ksp_type": "cg",
        "ksp_rtol": 1e-3,
        "pc_type": "hypre",
        },
        },
% For each concentration block k, using the ASM preconditioner:
"fieldsplit_[k-1]": {
    "ksp_rtol": 1E-1,
    "ksp_type": "gmres",
    "pc_type": "python",
    "pc_python_type": "firedrake.AssembledPC",
    "assembled": {
        "pc_type": "asm",
        "pc_asm_overlap": 1,
        "sub_pc_type": "ilu"
    },
    },

\end{lstlisting}

The reduction factor \texttt{REDFACTOR} for the coarse grid redistribution is
chosen depending on the number of processors such that the solver
sub-communicator always includes exactly 8 nodes.
The default overlap for ASM is one element.

\rme{
The custom preconditioner \texttt{CoarsenPenaltyPMGPC} is $p$-multigrid with $p=1$ on the coarse level where the proper penalization parameters are used for the DG discretization.
The custom preconditioner \texttt{CellIntegralPC} is used to ignore inter-element contributions.
These custom solvers are available in the \texttt{EchemFEM}~package~\cite{echemfem}.
}

\end{document}